\documentclass[12pt]{article}

\usepackage{amsmath,amsfonts,amssymb,mathrsfs,array}
\usepackage{graphicx,color}
\usepackage{theorem}
\newtheorem{thm}{Theorem}[section]
\newtheorem{lem}[thm]{Lemma}

\newtheorem{cor}[thm]{Corollary}

\newtheorem{prop}[thm]{Proposition}
\theorembodyfont{\upshape}
\newtheorem{proof}{Proof}

\newtheorem{defn}[thm]{Definition}
\newcommand{\qed}{\nopagebreak\par\hspace*{\fill}$\square$}

\newcommand{\ZZ}{{\mathbb Z}}

\newcommand{\CC}{{\mathbb C}}
\newcommand{\NN}{{\mathbb N}}
\def\gg{{\mathfrak{g}}}

\def\hh{{\mathfrak{h}}}

\def\LL{{\mathfrak{L}}}
\numberwithin{equation}{section}

\title{Integrable modules for Lie tori}
\author{S.Eswara Rao, Sachin S. Sharma}
\date{}

\begin{document}

\maketitle

\begin{abstract}
 In this paper we consider the universal central extension of a centerless Lie torus and 
 classify its irreducible  integrable 
modules  when the center acts non-trivially. They turn out to be highest weight modules for the direct 
sum of finitely many affine Lie algebras upto an automorphism.
\end{abstract}

\section*{Introduction}
Centerless Lie tori play an important role in explicitly constructing the extended affine Lie algebras(EALA);
they play similar role as derived algebras modulo center play in the realization of affine Kac-Moody
algebras \cite{Y, EN3}. Yoshii and Neher \cite{YIYJ, EN2} studied Lie torus in detail mainly due 
to its connection with
EALA. Starting with a centerless Lie torus, Erahard Neher \cite{EN3} has completed realization 
of EALA. In \cite{JS} using descent construction J. Sun obtained a universal central extension of Lie torus. 

In this work we use the multiloop realization of Lie tori given in
terms of multiloop algebras \cite{ABFP1}. 
These algebras cover almost all Lie tori except 
one class of Lie tori which appear in type A and coming from quantum tori \cite{ABFP1}.
We start with a centerless Lie torus and 
consider its universal central extension denoted by $\overline{LT}$. 
We add a finite set derivations of zero degree which 
measures the natural gradation on $\overline{LT}$, and denote it by $\stackrel{\sim}{LT}$. 

The purpose of the paper is to classify the irreducible integrable modules 
for $\stackrel{\sim}{LT}$ with finite dimensional 
weight spaces when the center acts non-trivially. 
They turn out to be highest weight modules for direct sum of 
finitely many affine Kac-Moody Lie algebras. Some very special 
cases are done in \cite{E3, EB, EZ} and \cite{XT}. 

 Let $LT$ be the centerless Lie torus .
Let $\overline{LT}$ be the universal central extension of $LT$ (Proposition \ref{prop1}). 
Both  $LT$ and  $\overline{LT}$ are 
naturally $\ZZ^{n+1}$ graded for some $n$. Let $\overline{D}$ 
be the spaces spanned by zero degree derivations 
$d_0, d_1, \cdots, d_n.$ Let $\stackrel{\sim}{LT} = \overline{LT} \oplus \overline{D}.$ 
The role of the $\overline{D}$ is to keep 
track of the gradation for $\overline{LT}$.

The Lie torus is defined as multiloop algebra using finitely many commuting
finite order automorphisms (subsection \ref{subs2}). In the case where 
all automorphisms are trivial we get the standard non-twisted multiloop algebra. 
In this case $\stackrel{\sim}{LT}$ is nothing 
but the toroidal Lie algebra (subsection \ref{subs1}). The irreducible integrable modules for 
toroidal case has been classified in \cite{E3}.	

The classification of the irreducible integrable modules (where the zero degree 
center acts non trivially) for $\stackrel{\sim}{LT}$ runs parallel to the toroidal case. 
But there are several places where we need completely different arguments. 
For example, in the toroidal case, highest weight exists at a graded level;
but in the case of $\stackrel{\sim}{LT}$, we 
could only able to prove that the highest weight exists at a non-graded level.
Because of the presence of automorphisms, the 
proof becomes more  complex.

Let the zero degree center  $K_0, K_1, \cdots, K_n$ act as $C_0, C_1, \cdots, C_n,$ where
$C_i \in \ZZ$ for $0 \leq i \leq n$ (see \cite{E3}, Lemma 2.3).
Upto change of co-ordinates (subsection \ref{subs4}) we can assume 
$C_0 > 0$ and $C_i = 0, 1 \leq i \leq n.$
Let $\LL$ denote a subquotient of $\stackrel{\sim}{LT}$ which 
is obtained by removing certain derivations and going 
modulo a part of the center.  The modules for $\stackrel{\sim}{LT}$ and $\LL$ are
called graded and non-graded respectively (see sections \ref{sec7} and \ref{sec8}).

The paper is organised as follows. In section \ref{sec1} and \ref{sec2} we 
give definitions and properties of $LT$ and $\stackrel{\sim}{LT}$. In section
\ref{sec3} we defined highest weight irreducible modules for $\LL$ which need 
not have finite dimensional weight spaces, and give 
necessary and sufficient condition for them to have finite dimensional weight spaces 
(Theorem \ref{thm1}). Sections \ref{sec4}, \ref{sec5}, and \ref{sec6} are devoted
to prove that an irreducible integrable $\LL$-module, 
where the center $K$ acts by a positive integer is a 
highest weight $\LL$-module (Theorem \ref{thn2}), and an 
irreducible integrable highest weight $\LL$-module is actually a module 
for the direct sum of finitely many affine Lie algebras 
(Proposition \ref{pr2} and Lemma \ref{lem7}). Finally, in sections \ref{sec7} and \ref{sec8} 
we complete the study
of the irreducible integrable $\stackrel{\sim}{LT}$-modules by establishing 
a one to one correspondence between irreducible modules of $\stackrel{\sim}{LT}$ 
and $\LL$.
\section{Notation and Preliminaries}\label{sec1} 
All vector spaces, algebras and tensor products are over 
complex numbers $\CC$. Let $\ZZ$, $\NN$ and $\ZZ_+$ 
denote integers, non-negative integers and positive integers respectively.
\subsection{Multiloop algebras and Lie torus}\label{subs1} Let $\gg$ be a 
simple finite dimensional Lie algebra and let $(\cdot|\cdot)$ be a non-degenerate symmetric bilinear form 
on $\gg$. Fix a positive integer $n$ and let $\sigma_0, \sigma_1, 
\cdots, \sigma_n$ be commuting finite order automorphisms 
of $\gg$ of order $m_0, m_1, \cdots, m_n$ respectively. 
Let $m=(m_1, \cdots, m_n) $, $k=(k_1, \cdots, k_n)$ and 
$l=(l_1, \cdots, l_n)$ are vectors in $\ZZ^n$. 
Let $\Gamma = m_1 \ZZ \oplus \cdots \oplus m_n \ZZ$ and 
$\Gamma_0= m_0 \ZZ$. Let $\Lambda = \ZZ^n/\Gamma$ and $\Lambda_0
= \ZZ/\Gamma_0$. Let $\overline{k}, \overline{l}$ 
denote the images in $\Lambda$. For any integers $k_0$ and $l_0$, 
let $\overline{k}_0$ and $\overline{l}_0$  denote images 
in $\Lambda_0$.
Let 
$$
\begin{array}{llll}
A&=& \CC[t_0^{\pm1}, \cdots, t_n^{\pm1}],\\
A_n&=& \CC[t_1^{\pm1}, \cdots, t_n^{\pm1}],\\
A(m)&=& \CC[t_1^{\pm m_1}, \cdots, t_n^{\pm m_n}], \ \mbox{and}\\
A(m_0, m)&=& \CC[t_0^{\pm m_0}, \cdots, t_n^{\pm m_n}]
\end{array}
$$
be Laurent polynomial algebras with respective variables.

Let  $\gg \otimes A$ be a loop algebra with usual bracket operations.
For $k \in \ZZ^n$, let $t^k = t_1^{k_1} \cdots t_n^{k_n} \in A_n.$
Let $X (k_0, k) = X \otimes t_0^{k_0} t^k, X \in \gg, k_0 \in \ZZ, k \in \ZZ^n$
denote a typical element of $\gg \otimes A$.
We now recall the definition of toroidal Lie algebra.
  
 Let $\Omega_A$ be the vector space spanned 
by symbols $t_0^{k_0}t^k K_i, 0\leq i \leq n, k_0 \in \ZZ$ and $k \in \ZZ^n$.
Let $dA$ be the subspace of $\Omega_A$ spanned 
by $\displaystyle{\sum^n_{i=0}} k_i t_0^{k_0} t^k K_i$. 
Similar spaces can be defined for all other Laurent 
polynomial algebras.
Let $Z= \Omega_A/dA, Z(n) =  \Omega_{A_n}/dA_n, \ Z(m) = \Omega_{A(m)}/dA(m)$ and 
$Z(m_0, m) = \Omega_{A(m_0,m)}/d A(m_0, m).$ Notice that $Z(n), Z(m), Z(m_0, m)$ 
are all subspaces of $Z$.
Then  $\tau = \gg \otimes A \oplus \Omega_A/dA$ is 
called toroidal Lie algebra with the following bracket operations:
\begin{enumerate}
 \item[{(a)}] $[X(k_0,k), Y(l_0, l)] = [X, Y] (k_0+l_0, k+l) + (X|Y) 
\displaystyle{\sum^n_{i=0}} k_i t_0^{k_0+l_0} t^{k+l}K_i;$
 \item[{(b)}] $\Omega_A/dA$ is central.
\end{enumerate}
It is well known that $\tau$ is the universal central extension of $\gg \otimes A$.

Notation : Let $\gg_{1}$ be any Lie algebra and $\hh_{1}$ be its 
finite dimensional ad-diagonalizable subalgebra.
We set for $\alpha \in \hh_{1}^*,$ 
$$
\gg_{1, \alpha} = \{x \in \gg_1 \,|\, [h,x] = \alpha(h) x, h \in \hh_{1}^* \} .
$$
Then we have $$\gg_1 = \displaystyle{\bigoplus_{\alpha \in \hh_{1}^*}} \gg_{1, \alpha} .$$
Let $\Delta(\gg_1, \hh_1) = \{\alpha \in  \hh_{1}^* \,|\, \gg_{1, \alpha} \neq 0 \}$ which includes $0.$
Let $\Delta^\times(\gg_1, \hh_1) = \Delta (\gg_1, \hh_1)\backslash{0}$.

\subsection{}\label{subs2}
 We will now define multiloop algebra as a subalgebra of $\gg \otimes A$. For $0\leq i \leq n$, let $\xi_i$ 
denote a $m_i$-th primitive root of unity.
Let
$$
\gg(\overline{k}_0,\overline{k}) = \{X \in \gg \,|\, \sigma_i X = \xi_i^{k_i} X, 0 \leq i \leq n \}.
$$
Then 
$\displaystyle{\bigoplus_{{(k_0,k)} \in \ZZ^{n+1}}} \gg(\overline{k}_0,
\overline{k}) t_0^{k_0} t^k$ is called multiloop algebra. The finite dimensional irreducible 
modules of the multiloop algebra  are classified by Michael Lau \cite{ML}.

 Throughout this paper we assume that $\gg(\overline{0}, \overline{0})$ is a simple Lie algebra. 
We can choose a Cartan subalgebra $\hh$ of $\gg$ such that 
$\hh(0) \subseteq \hh$ where $\hh(0)$ is a Cartan subalgebra of
$\gg(\overline{0}, \overline{0})$. We always fix such a choice of $\hh(0)$ 
and $\hh$ throughout this paper. It is well known that 
$\Delta^\times_0 = \Delta^\times_0 (\gg{(\overline{0}, \overline{0})}, \hh(0))$ 
is an irreducible reduced finite root system and 
has at most two root lengths. Let $\Delta^\times_{0,sh}$ be the set of non-zero short roots.
Define

$
\Delta^\times_{0,en}= \begin{cases} 
\Delta^\times_0 \cup 2 \Delta^\times_{0, sh} \ \mbox{if} \  \Delta^\times_0 \ \mbox{of \ type}\,\, B_l\\
\Delta^\times_0 \  \mbox{is otherwise}
\end{cases}
$                                              
and \,\,\,$\Delta_{0,en}= \Delta^\times_{0,en} \cup \{0\}.$
\begin{defn}
An irreducible finite dimensional highest weight module for 
$\gg(\overline{0}, \overline{0})$ is said to have property 
$(M)$ if the highest weight is one of the following
\begin{enumerate}
 \item[{(a)}] highest root.
 \item[{(b)}] highest short root if $\Delta^\times_0$ is not simply laced.
 \item[{(c)}] twice highest short root if $\Delta^\times_0$ is of type $B_l$.
\end{enumerate}
\end{defn}
In all these three cases the weight system of the module is contained in $\Delta_{0,en}$. Further the multiplicity 
of a non-zero weight is one.
\subsection{} We now define Lie torus (centerless) which is the main object of our study.
\begin{defn}\label{def1} A multi-loop algebra 
$\displaystyle{\bigoplus_{(k_0, k) \in \ZZ^{n+1}}} \gg(\overline{k}_0,\overline{k}) t_0^{k_0} t^k$ is called  a 
Lie torus and denoted by $LT$ if 

\begin{enumerate}
 \item  $\gg(\overline{0}, \overline{0})$ is a simple Lie algebra
 \item  As  $\gg(\overline{0}, \overline{0})$ module, each 
$\gg(\overline{k}_0,\overline{k}) = U (\overline{k}_0,\overline{k}) \oplus V(\overline{k}_0,\overline {k})$ where 
$U (\overline{k}_0,\overline {k})$ is trivial module and either $V(\overline{k}_0,\overline {k})$ is zero or satisfy 
the property $(M)$.
\end{enumerate}
\end{defn}

 Properties of a Lie torus $LT$ :
\begin{enumerate}
\item[{(a)}]  $\Delta (\gg, \hh(0)) = \Delta_{0,en}$ if $\Delta_0$ is of type $B_l$ and 
$V(\overline{k}_0,\overline {k})$ is isomorphic to highest weight irreducible finite dimensional module with highest 
weight $2 \beta_s$ for some $(\overline{k}_0,\overline {k})\neq (\overline{0}, \overline{0})$. Here $\beta_s$ is the 
highest short root in $\Delta^\times_0$.
\item[{(b)}]  $\Delta (\gg, \hh(0)) = \Delta_0$ in all other cases.
\item[{(c)}]  For $\alpha \in \hh(0)^*$, let 
$\gg (\overline{k}_0,\overline {k}, \alpha) = \{g \in \gg(\overline{k}_0,
\overline {k}) \,|\, [h,g] = \alpha(h)g, h \in \hh(0)\}$
Then $\gg (\overline{k}_0,\overline {k}) = \displaystyle{\bigoplus_{\alpha \in \hh^*(0)}}
\gg (\overline{k}_0,\overline {k}, \alpha) $.
The dimension of $\gg (\overline{k}_0,\overline {k}, \alpha) \leq 1$ for $\alpha \neq 0$.
\item[{(d)}] For $\alpha \neq 0,\, \gg (\overline{k}_0,\overline {k}, \alpha),
\,\gg (-\overline{k}_0,-\overline {k}, -\alpha) $ and 
$\Big[\gg (\overline{k}_0,\overline {k}, \alpha), 
\gg (-\overline{k}_0,-\overline {k}, -\alpha)\Big] $ forms a $\mathfrak{sl}_2$ copy.
\end{enumerate}
For details see Proposition 3.2.5 of \cite{ABFP1}.
\subsection{}\label{subs4} Change of co-ordinates: We work with $n$ variables for notational convenience.
Let $G=GL(n,\ZZ)$ be the group 
of $n\times n$ matrices with entries in $\ZZ$ and determinant $\pm1$. 
Then the group $G$ acts  naturally on $\ZZ^n$.
Denote the action by $B.k$ for $B \in G$ and $k \in \ZZ^n$. Let $s_i=t^{B.e_i}$,
where $\{e_i\}$ is the standard $\ZZ$ basis of 
$\ZZ^n$. Then clearly $A_n= \CC[s^{\pm1}_{1}, \cdots, s^{\pm1}_{n}].$ 
Consider a map $B : \gg \otimes A_n \oplus Z(n) \rightarrow \gg \otimes A_n \oplus Z(n)$ given by 
$$
\begin{array}{lll}
B.X\otimes t^k &=& X \otimes t^{B.k}\\
B.d(t^k)t^s&=& d(t^{B.k})t^{B.s}
\end{array}
$$
for $X \in \gg$ and $k, s \in \ZZ^n$, and where $d(t^k)t^s = \sum k_i t^{k+s} K_i$. 
Then it is easy to see that $B$ defines an 
automorphism.
Suppose $B=(b_{ij})$ then one can verify that $B(t^k K_i)= \displaystyle{\sum_j} b_{ij} t^{B.k} K_j$.
In particular $B(K_i)= \displaystyle{\sum_j} b_{ij} K_j$.
Now define a new Lie torus $\overline{LT}(B)$ by replacing variables $t_i$ by $s_i$ 
and note that $B$ takes the Lie torus 
$\overline{LT}$ to $\overline{LT}(B)$.
This is what we call change of co-ordinates. We will use this change of 
co-ordinates in the paper without any mention and just say that 
``upto a choice of co-ordinates.''

\section{Universal central extension of Lie torus}\label{sec2}
Let $LT= \displaystyle{\bigoplus_{(k_0,k) \in \ZZ^{n+1}}} 
\gg(\overline{k}_0,\overline {k}) t_0^{k_0} t^k$ 
be a Lie torus as defined in section 1.
 Let $\overline{LT} = LT \oplus Z(m_0,m)$.
Define a Lie algebra structure on $\overline{LT}$ by 
\begin{enumerate}
 \item[{(a)}]  $[X(k_0,k), Y(l_0,l)] = [X,Y] (k_0+l_0, k+l) + (X | Y)
 \displaystyle{\sum_{i=0}^n} k_i t_0^{k_0 + l_0} t^{k+l} K_i;$
 \item[{(b)}]  $Z(m_0,m)$ is central in $\overline{LT}$.
\end{enumerate}
Notice that $(X|Y) \neq 0 \Rightarrow k+l \in \Gamma$ and $k_0+l_0 \in \Gamma_0$.
This follows from the standard fact that 
$(\cdot|\cdot)$ is invariant under $\sigma_i, 0 \leq i \leq n.$ 
This proves that the above Lie bracket is closed. Notice also the above Lie bracket
is restriction of the bracket defined in subsection \ref{subs1}.
\begin{prop}\label{prop1}
$\overline{LT}$ is the universal central extension of $LT$.
\end{prop}
\begin{proof} 
 See Corollary (3.27) of \cite{JS}.
\end{proof}
\subsection{}\label{subs3}
Both $LT$ and $\overline{LT}$ are naturally $\ZZ^{n+1}$ graded.
To reflect this fact we add derivations. Let $\overline{D}$ 
be the space spanned by $d_0, d_1, \cdots, d_n$.
 Let $\stackrel{\sim}{LT} = \overline{LT} \oplus \overline{D}$. 
 Extend the Lie bracket in the following way:
$$
\begin{array}{lll}
\text{[}d_i, X(k_0, k)\text{]} &=& k_i X (k_0, k);\\
\text{[}d_i, t_0^{k_0} t^k K_j\text{]} &=& k_i t_0^{k_0} t^k K_j;\\
\text{[}d_i, d_j\text{]} &=& 0.
\end{array}
$$
Notice that $Z(m_0,m)$ is no more central in $\stackrel{\sim}{LT}$
but only an abelian ideal. In fact any graded subspace of 
$Z(m_0,m)$ is an ideal.

 Let $\stackrel{\sim}{\hh}= \hh(0) \oplus \displaystyle{\sum_{0 \leq i
 \leq n}}\CC K_i \oplus \overline{D}$ 
which is an abelian subalgebra of  $\stackrel{\sim}{LT}$. Let $\delta_i \in \stackrel{\sim}{\hh}^*$ such that 
$\delta_i (\hh(0)) = \delta_i (K_j) = 0$ and $\delta_i(d_j) = \delta_{ij}, 0 \leq i, j \leq n.$ For 
$(k_0, k) \in \ZZ^{n+1}$, let $ \delta_k = \sum k_i \delta_i$ and $\delta(k_0, k) = k_0 \delta_0 + \delta_k$.
Let
$$
\stackrel{\sim}{LT}_{\alpha + \delta(k_0,k)}=
\begin{cases}
\gg({\overline{k}_0, \overline{k}, \alpha)} \otimes t^{k_0}_0 t^k, \alpha \neq 0,\\
\gg({\overline{k}_0, \overline{k}, 0)} \otimes t^{k_0}_0 t^k \oplus 
\displaystyle{\sum_{\substack {\overline{k}_0, \overline{k} 
\in \Gamma_0 \oplus \Gamma\\{0 \leq i \leq n}}}} { t^{k_0} t^k K_i}, \,\,\alpha = 0
\,\, \mbox{and}\,\,  (k_0, k) \neq (0,0),\\
\stackrel{\sim}{\hh} = \gg(\overline{0}, \overline{0}, 0) \oplus 
\displaystyle{\sum_{0 \leq i \leq n}} \CC K_i \oplus 
\overline{D},\,\, (k_0, k, \alpha) = (0,0,0).\\
\end{cases}
$$
Then $\stackrel{\sim}{LT} = 
\displaystyle{\bigoplus_{\substack{\alpha \in \hh(0)^*\\(k_0, k) 
\in \ZZ^{n+1}}}} \stackrel{\sim}{LT}_{\alpha + \delta (k_0, k)}$ 
is a root space decomposition with respect to $\stackrel{\sim}{\hh}$ 
and each root space is finite dimensional.
\begin{defn} 
 A module $V$ of $\stackrel{\sim}{LT}$ is called weight module if 
\begin{enumerate}
 \item[{(a)}]  $V= \displaystyle{\bigoplus_{\lambda \in \stackrel{\sim}{\hh^*}}} V_\lambda,  
V_\lambda = \{v \in V \,|\, h v = \lambda (h)v, h \in \stackrel{\sim}{\hh} \}.$
 \item[{(b)}]  $\dim V_\lambda < \infty.$
\end{enumerate}
\end{defn}
\begin{defn}
A root $\alpha + \delta(k_0, k)$ of $\stackrel{\sim}{LT}$ is called real root if $\alpha \neq 0$ 
and null root if $\alpha = 0.$
\end{defn}
\begin{defn} \label{def2}
 A weight module $V$ of $\stackrel{\sim}{LT}$ is called integrable if every real root vector acts locally 
nilpotently  on $V$, i.e., for $X \in \stackrel{\sim}{LT}_{\alpha + \delta(k_0, k)}, 
\alpha \neq 0$ and $v \in V$, there exists 
$b=b(v, \alpha + \delta(k_0,k))$ such that $X^b. v =0$.
\end{defn}
The purpose of this paper is to classify irreducible integrable modules for $\stackrel{\sim}{LT}$. 
We will first reduce the 
problem to a subquotient of $\stackrel{\sim}{LT}$ and then classify those modules.

\section{Subquotient of $\stackrel{\sim}{LT}$ and highest weight module}\label{sec3} 
\subsection{} We will first define a Lie-algebra $\LL = LT \oplus \CC K \otimes A(m) \oplus \CC d_0$ 
where $K$ is a symbol. Let 
$X(k_0,k) \in \gg (\overline{k}_0, \overline{k})$ and
$Y(l_0, l) \in  \gg(\overline{l}_0, \overline{l})$ and define the bracket operations on $\LL$ as follows:
\begin{enumerate}
 \item[{(a)}]  $[X(k_0, k), Y(l_0, l)] = [X,Y] (k_0 + l_0, k+l) + 
 (X|Y)\, \delta_{l_0+k_0,0}\, k_{0} \,K \otimes t^{l+k};$
 \item[{(b)}]  $K \otimes A (m)$ is central;
 \item[{(c)}] $[d_0, X(k_0,k)] = k_0 X(k_0,k)$.
\end{enumerate}
 Let $\stackrel{\sim}{\LL} = \LL \oplus D$ 
where $D$ is the space spanned by derivations 
$d_1, d_2, \cdots , d_n$. Extend the Lie bracket 
to $\stackrel{\sim}{\LL}$ by defining $D$ action on $\LL$ as in subsection \ref{subs3}. We will 
first give a Lie algebra surjective homomorphism 
from $\Phi : \stackrel{\sim}{LT} \rightarrow \stackrel{\sim}{\LL}$ by 
$$
\begin{array}{lll}
 \Phi X(k_0,k) &=& X(k_0,k), X \in \gg (k_0, k);\\
 \Phi(t^{k_0}_0 t^k K_i) &=& 0 \ \mbox{if} \ i\neq 0 \ \mbox{or} \ k_0 \neq 0;\\
 \Phi(t^k K_0) &=& K \otimes t^k;\\
 \Phi(d_i) &=& d_i,0 \leq i \leq n.
\end{array}
$$
It is easy to see that $\Phi$ is a Lie algebra homomorphism. 
In fact $\ker \Phi$ is a graded subspace of $\ZZ(m_0,m)$ 
hence an ideal in $\stackrel{\sim}{LT}$ as noted in subsection \ref{subs3}.

We will now indicate the plan of the rest of the paper. 
The first step is to prove that any irreducible weight 
module $V$ of $\stackrel{\sim}{LT}$ is actually an irreducible 
module for $\stackrel{\sim}{\LL}$ upto a suitable choice of 
co-ordinates. The second step is to prove that any 
irreducible weight module for $\stackrel{\sim}{\LL}$  give 
rise to an irreducible module for $\LL$ with finite dimensional weight 
spaces with respect to 
$\stackrel{\sim}{\hh}\!\!(0) = \hh(0) \oplus \CC K \oplus \CC d_0$. 
We will also indicate how to recover the original module 
for $\stackrel{\sim}{LT}$ from the irreducible module of $\LL$. 
The third step is to prove that an irreducible 
integrable module for $\LL$ where $m_0 K$ acts as positive integer 
is an highest weight module (definition given below).
The last step is to classify such highest weight integrable irreducible modules.
\subsection{}
In this remaining part of section we will first take the last step.
 Let $\delta \in \stackrel{\sim}{\hh}\!\!(0)^*$ such 
that $\delta (\hh(0)) = 0, \delta(K)=0$ and $\delta(d_0) =1.$
Let
$$
\LL_{\alpha+k_0 \delta} =
\begin{cases}
\displaystyle{\bigoplus_{k \in \ZZ^n}} \gg (k_0, k, \alpha) t_0^{k_0} 
t^k, \ \mbox{if} \ \alpha+k_0 \delta \neq 0,\\\\
\displaystyle{\bigoplus_{k \in \ZZ^n}} \gg (0, \overline{k}, 0) 
t^k \oplus K \otimes A(m) \oplus \CC d_0,
\ \mbox{if} \ \alpha+k_0 \delta= 0.
\end{cases}
$$
Then $\LL = \displaystyle{\bigoplus_{\substack{\alpha \in \hh(0)^*\\k_0 \in \ZZ}}}
\LL_{\alpha+k_0 \delta}$ is  a root 
space decomposition with respect to $\stackrel{\sim}{\hh}(0)$.
But the dimensions of the root spaces are infinite dimensional.
 Recall that $\Delta_0 = \Delta (\gg(\overline{0},
\overline{0}), \hh(0))$ and $\Delta(\gg, \hh(0))$ from section \ref{sec1}. 
Let $\alpha_1, \cdots, \alpha_d$ be a  set of simple 
roots in $\Delta_0$. Let $\beta_0$ be a maximal root in 
$\Delta(\gg, \hh(0)).$ Let $\alpha_0 = -\beta_0+\delta$ 
which may not be root of $\LL$. 
Let $Q = \displaystyle{\bigoplus^d_{i = 0}} \ZZ 
\alpha_i$ and define an order $\leq$ on $\stackrel{\sim}{\hh}\!\!(0)^*$ for 
$\lambda, \mu \in \stackrel{\sim}{\hh}\!\!(0)^*, \,\lambda \leq \mu$ if 
$\mu-\lambda \in Q^+ = \displaystyle{\bigoplus^d_{i = 0}}\, \NN \alpha_i$. 
Clearly this ordering is consistent with the standard 
order on $\Delta(\gg, \hh(0))$. Let $\stackrel{\sim}{\Delta}$ be the set of all roots of $\LL$.
$$
\begin{array}{lll}
\stackrel{\sim}{\Delta}_+ = \{\alpha+k_0 \delta
\in \stackrel{\sim}{\Delta} \,|\,  k_0 > 0 \ \mbox{or} \ k_0=0, \alpha > 0 \}\\
\stackrel{\sim}{\Delta}_- = \{\alpha+k_0 \delta \in 
\stackrel{\sim}{\Delta} \,|\, k_0 < 0 \ \mbox{or} \ k_0=0, \alpha < 0 \}.\\
\end{array}\\
$$
Clearly $\stackrel{\sim}{\Delta}_+ $ is the set of positive root with the above ordering.
 Let
$$
\begin{array}{lll}
\LL^+ &= & \displaystyle{\bigoplus_{\alpha + k_0 \delta > 0}} 
\LL_{\alpha + k_0 \delta} , \\
\LL^- &= & \displaystyle{\bigoplus_{\alpha + k_0 \delta < 0}} \LL_{\alpha + k_0 \delta}, \,\mbox{and}\\
\LL^0 &= & \LL_0.
\end{array}\\
$$
Then clearly $\LL = \LL^+ \oplus \LL^0 \oplus \LL^- $ and note that $\LL^0$ 
is a subalgebra but need not be abelian. The bracket 
in $\LL^0$ does not produce the center $K \otimes A(m).$ 
In other words $K \otimes A(m)$ is a direct summand of $\LL^0$ as 
Lie algebras.
\subsection{}
 Construction of highest weight module for $\LL:$
Let $N$ be an irreducible finite dimensional module for $\LL^0$. Since $\stackrel{\sim}{\hh}\!\!(0) 
+ K \otimes A(m)$ is central, 
it is easy to see they act by scalars on $N$. Thus $\stackrel{\sim}{\hh}\!\!(0)$ acts on $N$ by a 
single linear function. 
Let $U(\LL)$ denote the universal enveloping algebra of $\LL$. Define Verma module for $\LL$.
$$
M(N) = U(\LL) \displaystyle{\bigotimes_{\LL^+ \oplus \LL^0}} N
$$
where $\LL^+$ acts trivially on $N$. By standard arguments, one can see that $M(N)$ admits a unique 
irreducible quotient say $V(N)$.
It is easy to see that $M(N)$ and $V(N)$ are weight module with respect to $\stackrel{\sim}{\hh}\!\!(0)$.
But they may not have 
finite dimensional weight spaces.
 We will now give a necessary and sufficient condition for $V(N)$ 
to have finite dimensional weight spaces.

Let $I$ be an ideal in $A(m)$ and let 
$$
\begin{array}{lll}
\LL^0 (I) &= & \displaystyle{\bigoplus_{k \in \ZZ^{n}}} 
\gg (0, \overline{k}, 0) t^k I \,\,\,\,\oplus\,\, K \otimes I,\\\\
\LL (I) &= & \displaystyle{\bigoplus_{\substack{k_0 \in \ZZ \\ k \in \ZZ^{n}}}}
\gg (\overline{k}_0, \overline{k}) \otimes t^{k_0}_0 t^k I \,\,\,\,\oplus \,\,K \otimes I.\\
\end{array}\\
$$
Clearly $\LL^0 (I)$ is an ideal in $\LL^0$ and $\LL(I)$ is an ideal in $\LL$.
Define $$\LL_{k_0} = \displaystyle{\bigoplus_{k \in \ZZ^n}} 
\gg (\overline{k}_0, \overline{k}) t^{k_0}_0 t^k, \,\, 
\mbox{and}\,\,\,\, \LL_{-k_0, -\alpha} = \displaystyle{\bigoplus_{k \in \ZZ^n}} 
\gg (-\overline{k}_0, \overline{k},
-\alpha) t^{-k_0}_0 t^k,\, k_0 \in \ZZ  .$$
Define $\LL_{k_0} (I)$ and $\LL_{-k_0, -\alpha} (I)$ by 
replacing $t^k$ by $t^k I$ in $\LL_{k_0}$ and  $\LL_{-k_0, -\alpha}$ respectively.
\begin{prop}\label{pr1} 
Suppose $N$ is finite dimensional irreducible module for $\LL^0$ such that $\LL^0(I) . N=0$ 
for some 
ideal $I$ of $A(m)$. Then $\LL(I).V(N)=0$.
\end{prop}
\begin{proof}
 Suppose $\alpha + k_0 \delta \geq 0$, then $\alpha + k_0 \delta = \displaystyle{\sum^d_{i = 0}} 
 s_i \alpha_i, s_i \in \NN.$
Define  $\mathrm{ht}(\alpha + k_0 \delta) = \displaystyle{\sum^d_{i = 0}} s_i \geq 0.$\\
 {\bf claim:}\ $ \LL_{-k_0, -\alpha} (I) N = 0$ for all roots 
$\alpha + k_0 \delta \geq 0$.\\
We will prove the claim by induction on the $\mathrm{ht}(\alpha + k_0 \delta)$. The claim is trivial if 
$\mathrm{ht} (\alpha + k_0 \delta)=0$ 
as $\alpha + k_0 \delta = 0$ and $\LL_{-k_0, -\alpha} (I) N \subseteq \LL^{0}(I)N = 0. $
Assume the claim for all roots $\alpha + k_0 \delta$ 
such that $\mathrm{ht} (\alpha + k_0 \delta) \leq s$ for a fixed $s > 0.$ 
{\bf subclaim:}\  $\LL_{l_0,\stackrel{\sim}{\alpha}}. \LL_{-k_0, -\alpha} (I) . N = 0 \,\,
\forall \stackrel{\sim}{\alpha} + 
l_0 \delta > 0$.\\
Consider the root $\stackrel{\sim}{\alpha} + l_0 \delta > 0$ and 
$$\LL_{l_0,\stackrel{\sim}{\alpha}}. \LL_{-k_0, -\alpha} (I) . N
=\LL_{-k_0, -\alpha}(I). \LL_{l_0,\stackrel{\sim}{\alpha}} N + [\LL_{l_0,\stackrel{\sim}{\alpha}}, 
\LL_{-k_0, -\alpha}(I)]N $$
The first term of RHS is zero as $N$ is the highest weight space. There are four cases 
$\stackrel{\sim}{\alpha} -\alpha   + (l_0-k_0) \delta >0, <0, = 0$ or may not be a root 
for the second term. The second term is zero in the first case as $N$ 
is highest weight space. The second term is zero in the second case by induction. In the third case, 
$[\LL_{l_0,\stackrel{\sim}{\alpha}}, \LL_{-k_0,-\alpha}(I)] \subseteq \LL^0(I)$ and 
hence zero by assumption. In the fourth case the bracket 
itself is zero. This proves the subclaim.

From the subclaim it follows $\LL_{-k_0, -\alpha} (I) N$ generates a proper submodule 
inside an irreducible module $V(N)$ for 
$\alpha + k_0 \delta < 0.$ Hence $\LL_{-k_0, -\alpha} (I) . N =0$ for all $\alpha + k_0 \delta > 0$. 
The same is true for 
$\alpha + k_0 \delta < 0$ as $N$ is highest weight space. 
This proves $\LL(I).N=0$. Consider $S=\{ v \in V(N) \,|\, \LL(I).v =0\}$
which contains $N$. Since $\LL(I)$ is an ideal in $\LL$, one can check that $S$ is a $\LL$-module.
This proves $S=V(N)$ and 
hence $\LL(I).V(N)=0.$ This completes the proof of the proposition.
\qed
\end{proof}

\begin{thm}\label{thm1}
$V(N)$  has finite dimensional weight spaces with respect to $\stackrel{\sim}{\hh}\!\!\!\,(0)$ if and only if 
there exists a co-finite ideal $I$ of $A(m)$ such that $\LL^0(I).N=0$.
\end{thm}
\begin{proof}
Suppose there exists a co-finite ideal $I$ of $A(m)$ such that $\LL^0(I)N=0.$ 
Then by Proposition \ref{pr1}, $\LL(I) V(N)=0.$ Thus $V(N)$ is a module for $\LL (A/I) \oplus \CC d_0.$ 
Here $\LL(A/I)$ is a Lie algebra in an obvious 
sense and admits root space decomposition with respect to $\stackrel{\sim}{\hh}\!\!(0).$ 
Since $A/I$ is finite dimensional, each 
root space of $\LL (A/I)$ is finite dimensional. By standard arguments (using PBW Theorem) 
it follows that any highest weight 
module for $\LL (A/I)$ has finite dimensional weight spaces. In particular $V(N)$ has 
finite dimensional weight spaces.

Conversely, suppose $V(N)$ has finite dimensional weight spaces with respect to $\stackrel{\sim}{\hh}\!\!(0).$
Since $\gg$ is simple and hence perfect, it follows that
\begin{equation} \label{eq2}
 \gg = \displaystyle{\sum_{\substack{\overline{k}_0,\overline{l}_0 \in \Lambda_0 \\ \overline{k}, 
 \overline{l} \in \Lambda}}} [\gg(\overline{k}_0,\overline{k}),
\gg(\overline{l}_0, \overline{l})].
\end{equation}
It is now easy to see that $\LL^0(I)$ is contained in the span of 
$[\gg(\overline{k}_0,\overline{k}),\gg (\overline{l}_0, \overline{l})] t^{k+l} I$, where $k_0 + l_0 = 0$ and
$\overline{k}, \overline{l} \in \Lambda$ and $K \otimes I.$\\
{\bf Claim:}
There exists a polynomial $Q_i$ in $t^{m_i}_i, (1 \leq i \leq n)$ such that $K \otimes t^k Q_i = 0$ on $V(N)$ 
for all $k \in \Gamma.$

To see the proof of the claim note that the non-degenerate form $(,)$ on $\gg$ remains non-degenerate on 
$\gg(\overline{0},\overline{0})$. Further it remains non-degenerate on $\hh(0).$ This is a general fact for any simple 
Lie algebra. Thus there exists $h_1,h_2 \in \hh(0)$ such that $(h_1, h_2) \neq 0$. 
Fix a non-zero vector $v$ in $N$ and $k_0 > 0.$
Consider $\{h_2 t^{m_iq}_i t_0^{-m_0} v, q \in \NN \}$, which belongs to a single weight space of $V(N)$. 
Thus there exists a polynomial $Q_i$ in variable $t_i^{m_i}$ such that 
$
h_2 t_0^{-m_0} Q_i v = 0
$.
Consider for any $k \in \Gamma$, 
$$h_1 t_0^{m_0} t^k. h_2 t_0^{-m_0} Q_i v 
= h_2 t_0^{-m_0} Q_i  h_1 t_0^{m_0} t^k v + (h_1, h_2) m_0 K t^k Q_i v =0 .$$ 
Here and below we omit $\otimes$ for convenience.
The first term of RHS is zero as $v$ belongs to $N$ and $m_0 \delta$ is a positive root. Thus $K t^k Q_i v =0$ for all 
$k \in \Gamma$. This proves the claim.

Now let $B(N) = \{ v_1, v_2, \cdots, v_s \}$ be a basis of $N$. Fix $(k_0,k) \in \ZZ^{n+1}$ 
such that $0 < k_i \leq m_i$ for all 
$i$. Note that $k_0 > 0$ but $\overline{k}_0$ could be zero.
Let
$
B\big(\gg(-\overline{k}_0,\overline{k})\big) = \{X_1, X_2. \cdots, X_r \}
$
be a basis of $\gg(-\overline{k}_0,\overline{k})$ in such a way that they are all $\hh(0)$ weight vectors.
Consider, for a fixed $v_j, X_p$ and $1 \leq i \leq n$
$
\{X_p t_0^{-k_0} t^k t_i^{m_i q} v_j, q \in \NN \}
$
which belongs to a single weight space of $V(N),$ which is finite dimensional. Thus there exists a polynomial 
$S_i : = S_i (X_p, v_j, k, k_0)$ in variable $t_i^{m_i}$ such that $X_p t_0^{-k_0} t^k S_i v_j =0$.
Let $P_i = \Pi S_i (X_p, v_j, k_0, k) Q_i$, where the product runs over all $p,j,k_0$ and $k$
such that $ 1\leq p \leq r, 1\leq j \leq s, 0 < k_0 \leq m_0 \ \mbox{and} \  0 < k_i \leq m_i$.
Here $Q_i$ is a polynomial given in the above claim. Since the product is finite, $P_i$ is a polynomial in $t_i^{m_i}.$
Consider for any $l \in \ZZ^n$ and $v_j \in N$
\begin{align*}
\gg(\overline{k}_0,\overline{l}) t_0^{-k_0} t^l X_p t^{-k_0}_0 t^k P_i& v_j
= \\ 
& X_p t^{-k_0}_0 t^k P_i \,\gg(\overline{k}_0,\overline{l}) t_0^{k_0} t^l  v_j 
+ [\gg(\overline{k}_0,\overline{l}),X_p]\, t^{l+k} P_i v_j.
\end{align*}
The central term does not appear as it is zero by the claim. For this note that 
$Q_i$ is a factor of $P_i$. The LHS is zero as $S_i$ is factor of $P_i$.
The first term on the right hand side
is zero as $v_j$ belongs to $N$ and $k_0 \delta > 0$ .
Thus we have proved that
$$
\big[\gg(\overline{k}_0,\overline{l}),  \gg(-\overline{k}_0,\overline{k})\big] l^{l+k} P_i v_j =0
$$
for $l \in \ZZ^n$ and $(k_0, k) \in \ZZ^n$ such that $0 < k_i \leq m_i$ for all $i$.
We already know that $K t^k P_i=0$ for all $k \in \Gamma$ by the claim. Let $I$ be the 
ideal generated by ${ P_1, P_2, \cdots, P_n}$ 
inside $A(m)$. It is easy to see that $I$ is a co-finite ideal.
Now from equation \ref{eq2} it is easily follows that $\LL^0(I) . N=0.$ 
This completes the proof of theorem.
\end{proof}

\section{Integrable modules for $\LL.$}\label{sec4}

\begin{defn} A module $V$ of $\LL$ is called a weight module if 
$$
V= \displaystyle{\bigoplus_{\lambda \in \stackrel{\sim}{\hh}\!\!(0)^*}} 
V_\lambda, \,\,V_\lambda = \{v \in V \,|\, hv= \lambda (h) v, \forall \, h\in \stackrel{\sim}{\hh}\!\!(0)^*\}
\ \mbox{and} \ \dim V_\lambda < \infty .
$$
\end{defn}
Recall that the roots of $\LL$ are of the form $\alpha + k_0 \delta$ and is called real root if $\alpha \neq 0.$
\begin{defn} 
A weight module $V$ of $\LL$ is called integrable if all real root vectors are locally nilpotent on 
$V$. 
\end{defn}
See definition \ref{def2} for more details. In this section we will classify irreducible integrable highest weight 
module for $\LL$. 

Let $V(N)$ be an irreducible highest weight module for $\LL$. Throughout this section we assume 
$V(N)$ is integrable . In particular 
$V(N)$ has finite dimensional weight spaces with respect to $\stackrel{\sim}{\hh}\!\!(0)$.
Thus by Theorem \ref{thm1}, there exists a co-finite ideal $I$ of $A(m)$ such that $\LL^0(I).N=0$. 
We can assume that the ideal $I$ 
generated by polynomials $P_i$ in variable $t^{m_i}_i, 1 \leq i \leq n.$ We can 
further assume that the constant term is one.
Write $P_i(t^{m_i}_i) = \displaystyle{\prod_{j=1}^{q_i}} (t^{m_i}_i - a_{ij}^{m_i})^{b_{ij}}$
for some positive integers $b_{ij}$ and $q_i$.
Further $a_{ij}^{m_i} \neq a_{ij^\prime}^{m_i}$ for $j \neq j^\prime$.
Let $P_i^\prime(t_i^{m_i}) = \displaystyle{\prod_{j=1}^{q_i}} (t^{m_i}_i - a_{ij}^{m_i})$.
Let $I^\prime $ be the ideal of $A(m)$ generated by $P_i^\prime , 1 \leq i \leq n.$ 
Then clearly $I \subseteq I^\prime.$
\begin{prop}\label{pr2} 
 $\LL(I^\prime) V(N) = 0.$ 
\end{prop}

We need some Lemmas. In view of Proposition \ref{pr1} it is sufficient to prove that 
 $\LL^{0}(I^\prime) N = 0.$ Fix a positive root $\alpha + k_0 \delta$ such that $\alpha \neq 0$.
 Let $J$ be a subquotient of 
$A(m)$ of the form $J_1/J_2$ where $J_1$ is an ideal of $A(m)$ and $J_2$ is an ideal of $J_1$.
Recall $$\LL =\bigoplus_{\substack{\alpha \in {\hh}(0)^*\\k_0 \in \ZZ}}{\LL_{\alpha+k_0 \delta}}.$$
Fix a root $\alpha + k_0 \delta > 0.$ Let $\gg (J)$ be the linear space spanned by
\begin{enumerate}
\item[{(a)}] $\displaystyle{\bigoplus_{k \in \ZZ^n}} \gg(\overline{k}_0,\overline{k},\alpha) t_0^{k_0} t^k J,$ 
\item[{(b)}] $\displaystyle{\bigoplus_{k \in \ZZ^n}} \gg(-\overline{k}_0,\overline{k}, -\alpha) t_0^{-k_0} t^k J,$  
\item[{(c)}] $\displaystyle{\bigoplus_{k \in \ZZ^n}} \gg(2\overline{k}_0,\overline{k}, 2\alpha) t_0^{2k_0} t^k J,$ 
\item[{(d)}] $\displaystyle{\bigoplus_{k \in \ZZ^n}} \gg(-2\overline{k}_0,\overline{k}, -2\alpha) t_0^{-2k_0} t^k J,$ 
\item[{(e)}] $\Big[\displaystyle{\bigoplus_{k \in \ZZ^n}} \gg(\overline{k}_0,\overline{k},\alpha) t_0^{k_0} t^k J \,\,,
              \,\,\displaystyle{\bigoplus_{k \in \ZZ^n}} \gg(-\overline{k}_0,\overline{k}, -\alpha) t_0^{-k_0} t^k \Big],$ 
\item[{(f)}] $\Big[\displaystyle{\bigoplus_{k \in \ZZ^n}} \gg(2\overline{k}_0,\overline{k},2\alpha) t_0^{2k_0} t^k J \,\,,
              \,\,\displaystyle{\bigoplus_{k \in \ZZ^n}} \gg(-2\overline{k}_0,\overline{k}, -2\alpha) t_0^{-2k_0} t^k \Big].$ 
\end{enumerate}
Note that the second term in $(e)$ and $(f)$ do not have $J$. In case $2\alpha +2k_0 \delta$ is not a root, 
then the terms in $(c), (d)$ and $(f)$ do not occur. Now it is easy to check that $\gg(J)$ is a 
Lie algebra using Jacobi identity.

Let $\stackrel{\sim}{\gg}(J) = \gg(J) + \sum \gg(0,\overline{k}, 0) t^k J + K \otimes J \oplus \CC d_0$.
Note that $\LL^0(J) \subseteq \stackrel{\sim}{\gg}(J)$ and the terms in $(e)$ and $(f)$ 
are contained in  $\LL^0(J)$. 
It is direct checking that $\gg(J)$ is an ideal in $\stackrel{\sim}{\gg}(J)$.
Let $\stackrel{\sim}{N}$ be the $\stackrel{\sim}{\gg}(A(m))$-module 
generated by $N$. Note the following:
\begin{enumerate}
 \item[{(1)}]  $\stackrel{\sim}{N}$ is actually module for $\stackrel{\sim}{\gg} (A(m)/I)$, where
 $A(m)/I$ is finite dimensional.
 \item[{(2)}]  The terms in $(a)$ and $(c)$ in the definition of $\gg(A(m))$ acts trivially on $N$ as they 
 correspond to positive root spaces.
 \item[{(3)}]  $\LL^0 (A_m)$ leaves $N$ invariant. 
 \item[{(4)}]  The terms in $(b)$ and $(d)$ of $\gg(A(m))$ acts locally nilpotently on 
 $N$ as they correspond to real root spaces.
\end{enumerate}
Thus by $PBW$ theorem we conclude that $\stackrel{\sim}{N}$ is finite dimensional. Also note that 
$\stackrel{\sim}{N}$ is $d_0$ -invariant.

Let $\stackrel{\sim}{\gg} (A(m)) = \,\,\stackrel{\sim}{\gg} (A(m))_{-} \oplus \LL^0 (A(m))
\oplus \stackrel{\sim}{\gg}(A(m))_+$
be a natural decomposition. The positive part corresponds to 
$\alpha + k_0 \delta$ and $2\alpha + 2k_0 \delta.$ Similarly the 
negative part. So  $\stackrel{\sim}{N}$ is an highest weight module for
$\stackrel{\sim}{\gg}(A(m))$ and hence has a unique 
irreducible quotient $N_1$.
Let $\pi : \stackrel{\sim}{N} \rightarrow N_1$ be the natural map.
 Note that $\ker \pi \cap N = \{0\}$. 
\begin{lem}\label{lem1} 
Fix $k$ and $l \in \ZZ^n,$ and $k_0 \in \ZZ$. Let
$$
Y \in  \gg(-\overline{k}_0,\overline{l}, -\alpha) t_0^{-k_0} t^l,\,\,
X \in  \gg(\overline{k}_0,\overline{k}, \alpha) t_0^{k_0} t^k I^\prime,\,\,
H = [X,Y],  Y_1 = [H,Y].
$$
Suppose there exists a vector $v$ in $N_1$ such that $Xv = Y_1v = [Y_1,Y] v = 0$ 
and $Hv = \lambda v$ for some $\lambda \in \CC.$
Then 
\begin{enumerate}
 \item[{(1)}] $H Y^q v = \lambda Y^q v, q \geq 0.$ 
 \item[{(2)}] $X Y^q v =  q \lambda Y^{q-1} v , q > 0.$
\end{enumerate}
\end{lem}
\begin{proof} We have the following formula from \cite{K}
\begin{equation}\label{eq3}
Y^q Y_1 = \displaystyle{\sum_{s = 0}^q}\begin{pmatrix}q\\s\end{pmatrix} ((adY)^s Y_1 ) Y^{q-s}.
\end{equation}
We first note that $(adY)^s Y_1 = 0$ for $s \geq 2$ as $(s+1) \alpha$ is not a root.
We also have $Y^q Y_1 v =0.$
Thus formula \ref{eq3} becomes
$$
Y_1Y^q v + q[Y , Y_1] Y^{q-1} v = 0.
$$
But $[Y,Y_1] Y^{q-1} v$ is zero as $[Y,Y_1]$ commutes with $Y$ and kills $v$.
 This proves $[H,Y] Y^q v = 0$ for $q \geq 0.$ We will now prove (1) by induction on $q$. 
By assumption in the 
Lemma, it is true for $q = 0.$ Assume for $q \geq 1$ and consider
$$
\begin{array}{lll}
H Y^{q+1} v &=& H YY^q v = YHY^q v + [H,Y] Y^q v\\
&=& \lambda Y^{q+1} v. 
\end{array}
$$
This completes the induction step. We will now prove (2) again by induction on $q$. (2) is true for $q=1$ as 
$XYv = YXv + [X,Y] v = Hv = \lambda v$.
Assume for $q>1$ and consider
$$
\begin{array}{lll}
XY^{q+1} v &=& YXY^q v + HY^q v\\
&=& q \lambda Y^q v + \lambda Y^q v \ \mbox{(by 1 and induction)}\\
&=& (q+1) \lambda Y^q v.
\end{array}
$$
This completes the induction step. Thus the proof of the Lemma is completed.
\qed
\end{proof}
\begin{lem} 
$\stackrel{\sim}{\gg}(I^\prime / I)$ is a solvable Lie-algebra.
 \end{lem}
 \begin{proof}
 It is easy checking. Just note that $I^{\prime P} \subseteq I$ for large $P$.
 \qed
 \end{proof}
 \begin{lem}\label{lem8} 
 $\gg(I^\prime / I) N_1 = 0$.
  \end{lem}

\begin{proof}
Since $\gg(I^\prime / I)$ is solvable and $N_1$ is finite dimensional, there exists a vector $v$ in $V$ 
such that $\gg(I^\prime / I)$ acts as scalars.\\
{\bf Claim} \ $\gg(I^\prime / I) v = 0$.\\
The positive and negative part act as locally nilpotent operators and also as scalars. Hence they acts trivially.
It remains to prove the zero part of $\stackrel{\sim}{\gg} (I^\prime / I)$ acts trivially on $v$. 
Fix $k,l \in \ZZ^n$ and 
$X,Y,Y_1$ and $H$ be as in Lemma \ref{lem1}. We have already seen $Xv=0$, $Y_1v = [Y_1,Y] v =0$ as $Y_1$ and $[Y_1,Y]$ 
are in the negative part of  $\gg(I^\prime / I)$. We have $Hv=\lambda v$.
Thus we have Lemma \ref{lem1}.
We will prove that $\lambda = 0$. Let $P$ be least positive integer such that $Y^Pv=0$ and $Y^{P-1} v \neq 0$.
Consider
$$
\begin{array}{lll}
\lambda Y^{P-1} v &=& HY^{P-1}v\\
&=&[X,Y] Y^{P-1}v\\
&=&(XY-YX) Y^{P-1}v\\
&=& -(P-1) \lambda Y^{P-1}v .\\
\end{array}
$$
This means $P \lambda =0 \Rightarrow \lambda =0$ as $P \neq 0.$ Recall the definition of 
$\gg(I^\prime / I)$ and we have 
proved that the term from $(e)$ acts trivially on $v$. Similar argument proves that the 
term $(f)$ acts trivially. This 
completes the proof of claim.

 Now let $M=\{ v \in N_1 \,|\, \gg(I^\prime / I) v =0\}$. Since $\gg(I^\prime / I)$ is in ideal in 
$\stackrel{\sim}{\gg} (I^\prime / I), M$ will be $\stackrel{\sim}{\gg} (I^\prime / I)$-module. 
From above, we 
know that $M\neq 0.$ Since $N_1$ is $\stackrel{\sim}{\gg} (I^\prime / I)$ irreducible 
it follows that $N_1=M$ 
and the lemma follows.
\qed
\end{proof}
{\bf Proof of Proposition \ref{pr2}}\,\,  Note that as $\pi(N) \hookrightarrow N_1$, 
we have $\gg(I^\prime / I) N = 0$ by above lemma. 
So to prove
 $\LL^{0}(I^\prime)N=0$, we only need to show that 
$\gg(0,\overline{k},0) t^k I^\prime \oplus K \otimes I^\prime$ is zero on $N$ for all $k \in \ZZ^n$. 
From the above lemma the term
\begin{align*} \label{eq1}
\Big[&\gg(\overline{k}_0,\overline{k}, \alpha)t_0^{k_0} t^k I^\prime , 
\gg(-\overline{k}_0,\overline{l}, -\alpha) t_0^{-k_0} t^l\Big]\\
&=\Big[\gg(\overline{k}_0,\overline{k}, \alpha), \gg(-\overline{k}_0,
\overline{l}, -\alpha)\Big] t^{k+l} I^\prime 
+(\gg(\overline{k}_0,\overline{k}, \alpha), 
\gg(-\overline{k}_0,\overline{l}, -\alpha)) k_0 K\otimes t^{k+l}I^\prime 
\end{align*}
is zero on $N$.

The above statement is true for any positive root $\alpha + k_0 \delta$.
Now replacing $k_0$ by $m_0 + k_0$  we see that first 
term of RHS of the above equation does not change as 
$\overline{m_0+k_0} = \overline{k}_0$. By subtracting the new equation from the above 
equation we get 
$K \otimes t^k I^{\prime} N=0$ for all $k \in \Lambda$, as
 $(\gg(\overline{k}_0,\overline{k}, \alpha),\gg(-\overline{k}_0,\overline{k}, -\alpha))$
is non-zero.
 To complete the proof of the proposition,
it is sufficient to see that 
$\gg(0, \overline{k^\prime},0)$ is spanned by $\Big[\gg(\overline{k}_0,\overline{k}, \alpha), 
\gg(-\overline{k}_0,\overline{l}, -\alpha)\Big]$. 
Now this follows from identity \ref{eq2} and the observation that 
$\gg = \displaystyle{\bigoplus_{\alpha \neq 0 }} (\gg_\alpha+[\gg_\alpha,\gg_{-\alpha}])$.

\section{Affine Lie algebras.}\label{sec5}
 In the last section we have seen that the irreducible integrable module $V(N)$ is actually a module for 
$\LL(A(m)/I^\prime )$. In this section we will describe $\LL(A(m)/I^\prime )$ and prove that it is isomorphic 
to direct sum of finitely many affine Lie algebras.
\subsection{}
Recall that $\sigma_0$ is an automorphism of order $m_0$ and $\xi_0$ is $m_0$-th primitive root of unity.
Let 
$$
\gg_{\overline{k}_0}= \{x \in \gg \,|\, \sigma_0 x = \xi_0^{k_0} x \}.
$$
Write $$\LL(\gg, \sigma_0)= \displaystyle{\bigoplus_{k_0 \in \ZZ}} 
\gg_{\overline{k}_0} \otimes t^{k_0} \oplus \CC K$$
which is known to be an affine Lie algebra.\\
{\bf Some notations:}
For $1 \leq i \leq n$, let $N_i$ be a positive integer and let $N=N_1 \cdots N_n$ be the product.
Let $\underline{a}_i = ({a}_{i1}, \cdots, a_{i_{N_i}})$ be non-zero complex numbers such that 
$a^{m_i}_{ij} \neq a^{m_i}_{ij^\prime}$ for $j \neq j^\prime$. Let $I= (i_1,\cdots, i_n)$ such
that $1 \leq i_j \leq N_j$. 
There are $N$ of them and let $I_1,\cdots, I_N$ be some order.
Write 
$$
a_I^q = a_{1i_1}^{q_1} \cdots a_{ni_n}^{q_n}, q=(q_1, \cdots , q_n) \in \ZZ^n.
$$
Fix $k=(k_1,\cdots, k_n) \in \ZZ^n$ such that $0 \leq k_i < m_i$. 
\begin{lem}\label{lem9} Consider the $N \times N$ matrix 
$$
B=( a^{l_1m_1+k_1}_{1i_1} \cdots  a^{l_nm_n+k_n}_{ni_n})_{\substack{1\leq l_i \leq N_i\\1\leq i_j 
\leq N_j\\1\leq i,j \leq n.\\}}
$$
Then $B$ is invertible.
\end{lem}
\begin{proof} Let $B_i = (a_{ij}^{m_il_i})_{\substack{1\leq j \leq N_i\\1\leq l_i \leq N_i}} \ \mbox{and}$
let $D_i = \begin{pmatrix} a_{i1}^{k_i}  &&  0 \\ &\ddots & \\0 && a_{iN_i}^{k_i} \end{pmatrix}$.
Since $B_i$ is Vandermonde matrix, it is invertible. $D_i$ is  also invertible.
Then it is not too difficult to see $B$ is a tensor product of $B_iD_i$. It is well known that the tensor product 
of invertible matrices is invertible; hence $B$ is invertible. See \cite{B} for more details. 
\qed
\end{proof}
Let $\LL^\prime = \displaystyle{\bigoplus_{(k_0,k) \in \ZZ^n}} \gg(\overline{k}_0,\overline{k})t_0^{k_0}t^k \oplus
K \otimes A(m) $ and define
$\LL^\prime (I)$ similarly.
Define a Lie algebra homomorphism
$\varphi : \LL^\prime \rightarrow \displaystyle{\bigoplus_{N-copies}} \LL (\gg, \sigma_0)$ given by
\begin{enumerate}
\item $\varphi(Xt_0^{k_0}t^k) = (X t^{k_0} a^k_{I_1}, \cdots,X t^{k_0} a^k_{I_N}).$
\item $\varphi(K t^{k}) = (Ka^k_{I_1}, \cdots, Ka^k_{I_N})$, where $k \in \Gamma$.
\end{enumerate}
Let $P_i^\prime(t_i^{m_i}) = \displaystyle{\prod_{j=1}^{N_i}} (t_i^{m_i}-a_{ij}^{m_i})$ and 
let $I^\prime$ be an ideal generated 
by ${ P_1^\prime,\cdots, P_n^\prime}$ inside $A(m)$.
\begin{lem}\label{lem7}
\begin{enumerate}
 \item[{(1)}] $\varphi$ is surjective.
 \item[{(2)}] $\LL^\prime/\LL^\prime (I^\prime) \cong \displaystyle{\bigoplus_{N-copies}} \LL (\gg, \sigma_0)$.
\item[{(3)}]  $\varphi$ is $d_0$ grade preserving map.
\end{enumerate}
\end{lem}
\begin{proof} Fix $k \in \ZZ^n$ such that $0 \leq k_i < m_i$. Fix any $k_0 \in \ZZ$.
Let $T(\overline{k}_0,\overline{k})$ be the span of $\gg(\overline{k}_0,\overline{k}) t_0^{k_0} t^k t^Q$ where 
$Q=(q_1m_1,\cdots, q_n m_n) \in \Gamma$ and $1 \leq q_i \leq N_i$. Then by Lemma \ref{lem9}, the restriction of 
$\varphi$ to 
$$
T(\overline{k}_0,\overline{k}) \rightarrow \oplus \gg(\overline{k}_0,\overline{k}) t^{k_0}
$$
is both injective and surjective as the corresponding matrix is invertible. Similar argument holds 
good for $\varphi$ 
restricted to $K \otimes A(m)$. This proves (1). To see (2), first note that 
$\LL^\prime(I^\prime) \subseteq \ker \varphi$. 
It is easy to see that $T(\overline{k}_0,\overline{k})$ is a spanning set for  
$\gg(\overline{k}_0,\overline{k}) t_0^{k_0} t^k A(m)/I^\prime$. In fact it is a basis as 
the map  $\varphi$ is injective on 
$T(\overline{k}_0,\overline{k})$. It proves $\LL^\prime(I^\prime)= \ker \varphi$.
(3) is obvious.
\qed
\end{proof}
 {\bf Remark:}
By a result of Kac \cite{K} we can assume that $\LL(\gg, \sigma_0)$ is an affine Lie algebra.
On the surface the highest weight modules that appear here look different from the standard highest
weight modules in \cite{K}. The reason is ad-diagonalizable subalgebra $\stackrel{\sim}{\hh}\!\!(0)$ is much smaller
than the Cartan subalgebra of $\LL (\gg, \sigma_0)$. But because of integrability it 
is easy to see that they are same. First we note 
that an irreducible highest weight module for finitely many copies of 
$\LL (\gg, \sigma_0)$ is actually tensor product of highest weight modules for 
$\LL (\gg, \sigma_0)$. Any highest module for $\LL (\gg, \sigma_0)$ is automatically decomposes under $d$ action;

Suppose we start with a standard highest weight irreducible integrable module for $\LL (\gg, \sigma_0)$; 
see \cite{K}. It is standard fact that each $d$-eigenspace is finite dimensional. 
 In particular it has finite dimensional weight spaces. 
Conversely, any highest weight module is a standard highest weight module as the 
top component is finite dimensional and so the 
weights are bounded above. 

\section{Integrable modules revisited}\label{sec6}
Recall the definition of integrable modules for $\LL$ from section \ref{sec4}. 
In this section we will use techniques in \cite{REC} to prove that an irreducible 
integrable module for $\LL$ where $m_0K$ acts as positive integer is an highest weight module.
\subsection{}
Let $\gg_{aff}= \gg(\overline{0},\overline{0}) \otimes 
\CC [t_0^{m_0},t_0^{-m_0}] \oplus \CC K \oplus \CC d_0$ which is a 
subalgebra of $\LL$. Recall that $\stackrel{\sim}{\hh}\!\!(0) = \hh(0) \oplus 
\CC K \oplus \CC d_0 \subseteq \gg_{aff}$.
\begin{thm}\label{thn2}  Suppose $V$ is an irreducible module for $\LL$ with 
finite dimensional weight spaces with respect to 
$\stackrel{\sim}{\hh}\!\!(0)$. Further, suppose $V$ is integrable for $\gg_{aff}$
where the canonical central element $m_0K$ acts as 
positive integer, then $V$ is an highest weight module for $\LL$.
\end{thm}
Before proving the above theorem, we need some notations and some lemmas.

 Recall that $\delta \in \stackrel{\sim}{\hh}\!\!(0)^*$ such 
 that $\delta (\hh(0))=0, \delta(K) = 0$ and 
$\delta(d_0)=1$. Now define $w \in \stackrel{\sim}{\hh}\!\!(0)^*$ 
such that $w (\hh(0))=0, w(K)=1$ and $w(d_0) = 0$.  
Given $\lambda \in \stackrel{\sim}{\hh}\!\!(0)^*$, let $\overline{\lambda}$ be the restriction
to $\hh(0)$. Conversely, given any $\overline{\lambda} \in {\hh}(0)^*$ 
we can extend $\overline{\lambda}$ to $\stackrel{\sim}{\hh}\!\!(0)$ such that 
$\overline{\lambda}(K) = \overline{\lambda} (d_0) = 0$. Given any $\lambda 
\in \stackrel{\sim}{\hh}\!\!(0)^*$, we can write in a unique
way 
$
\lambda = \overline{\lambda} + \lambda(d_0) \delta + \lambda (K) w.
$

Recall that $\Delta_0 = \Delta (\gg(\overline{0},\overline{0}), \hh(0))$ 
and $\Delta = \Delta (\gg, \hh(0))$. Let 
$\beta$ be maximal root in $\Delta_0$ and let $\beta_0$ be maximal
root in $\Delta$. We have already note that 
$\beta = \beta_0$ except possibly in the case $B_l$. In that case $\beta_0 = 2\beta s$ where $\beta s$ 
is a highest short 
in $\Delta_0$. In any case $\beta \leq_0 \beta_0$ where $\leq_0$ denote the standard ordering in $\Delta_0$.
 Let $\alpha_0 = -\beta_0 + \delta$ and $\alpha^ 1= -\beta +m_0 \delta$. Note that $\alpha_0$ may not be a
root of $\LL$.
Recall that $\alpha_1, \alpha_2, \cdots, \alpha_d$ is a set of simple roots in $\Delta_0$. Let 
$\overset{\circ}{Q} = \displaystyle{\bigoplus_{i=1}^d} \,\ZZ \alpha_i, \,
Q_{aff}=\overset{\circ}{Q} + \ZZ \alpha^1 $ and 
$Q = \overset{\circ}{Q} + \ZZ \alpha_0$.
Let $\overset{\circ}{Q}^+ = \displaystyle{\bigoplus_{i=1}^d} \,\NN \alpha_i, 
\,Q^+_{aff} = \overset{\circ}{Q}^+ + \NN \alpha^1$ and 
$Q^+ = \overset{\circ}{Q}^+ + \NN \alpha_0$ .
Consider $\alpha^1 = -\beta +m_0 \delta = -\beta_0 + \delta +(\beta_0 -\beta) + (m_0-1) (\beta_0 + \alpha_0)$.
From this we can see 
$$
\overset{\circ}{Q} \subseteq Q_{aff} \subseteq Q \ \mbox{and } \ \overset{\circ}{Q}^+ 
\subseteq Q^+_{aff} \subseteq Q^+ 
.$$
\subsection{}
Define three orderings on $\stackrel{\sim}{\hh}\!\!(0)^*$ in the following way. Let $\lambda, 
\mu \in \stackrel{\sim}{\hh}\!\!(0)^*$
$
\begin{array}{lll}
\lambda \leq_0 \mu &\mbox{if}& \mu-\lambda \in \overset{\circ}{Q}^+,\\
\lambda \leq \mu &\mbox{if}& \mu-\lambda \in Q^+_{aff},\\
\lambda \leq_1 \mu &\mbox{if}& \mu-\lambda \in  Q^+.\\
\end{array}
$\\
From above we have $\lambda \leq_0 \mu \Rightarrow \lambda \leq \mu$
and $\lambda \leq \mu \Rightarrow \lambda \leq_1 \mu.$
 We note the following useful observation:
 \begin{equation}\label{ob1}
 \lambda \geq_1 0 \,\, \mbox{and}\,\, \lambda(d_0) = 0 \Longrightarrow \lambda \geq_0 0.
\end{equation}
 Let $\Delta_{aff} = \{ \alpha +k_0 m_0 \delta \,| \, k_0 \in \ZZ, \alpha \in \Delta_{0}\}$.
Let
$$
\begin{array}{lll}
\Delta^+_{aff} &=&  \{ \alpha +k_0 m_0 \delta \in \Delta_{aff} \, | \, k_0 > 0 \ \mbox{or} \  k_0=0, \alpha >_0 \,0\}\\
&=&  \{ \alpha +k_0 m_0 \delta \in \Delta_{aff} \,|\, \alpha +k_0 m_0 \delta > 0\}.
\end{array}
$$
It is obvious that $\gg_{aff}$ has root space decomposition
$
\gg_{aff} = \displaystyle{\bigoplus_{\alpha \in \Delta_{aff}}}
\gg_{aff, \alpha}$ and  $\gg_{aff,0} = \stackrel{\sim}{\hh}\!\!(0)
$.
Let $N_{aff}^+ = \displaystyle{\bigoplus_{\alpha> 0}}\gg_{aff,\alpha} $.
\subsection{} Let $V(\mu)$ denote the irreducible integrable highest weight module 
for $\gg_{aff}$ with highest weight $\mu$. The following are well known facts:\\
{\bf Fact 1 :} An highest weight integrable module for $\gg_{aff}$ is necessarily irreducible.\\
{\bf Fact 2 :} Suppose $\lambda$ and $\mu$ are dominant integral weights for the simple roots 
 $\alpha^1, \alpha_1, \cdots, \alpha_d$ 
such that $\lambda \leq \mu$, then $\lambda$ is a weight of $V(\mu)$.

 Recall that $\overset{\circ}{Q}$ is a root lattice of $\gg(\overline{0},\overline{0})$. Let 
$\overset{\circ}{\Lambda}$ be the weight lattice of $\gg(\overline{0},\overline{0})$ 
and let ${\overset{\circ}{\Lambda}}^+$ 
be the set of dominant integral weights.
\begin{defn} $\nu \in {\overset{\circ}{\Lambda}}^+$ is called minimal if 
 $\mu \leq_0 \nu, \ \mu \in {\overset{\circ}{\Lambda}}^+$, then $\mu = \nu$.
\end{defn}
 It is known that every coset $\overset{\circ}{\Lambda}/\overset{\circ}{Q}$ has a unique minimal weight.
Let $P(V(\mu))$ denote the set of weight of $V(\mu)$ and let 
$\overline{P}(\mu) = \left\{\overline{\nu} \,|\, \nu \in P(V(\mu))\right\}$. 
Clearly $\overline{P}(\mu)$ 
determines a unique coset in $\overset{\circ}{\Lambda}/\overset{\circ}{Q}$ . 
Let $\overline{\nu}_0$ be the minimal weight and note that 
$\overline{\nu}_0 \leq_0 \overline{\mu}$. Let $s$ be any complex number 
such that $\lambda(d) - s$ is a non-negative integer divisible by $m_0$.
It follows that $\nu_0 = \overline{\nu}_0 + s \delta + \mu(K) w \in P (L(\mu))$; 
see \cite{E3} for details. 

We fix an irreducible integrable module $V$ for $\LL$.
We assume that $m_0K$ acts as positive integer. Let $P(V)$ be the set of 
weights. Let $P_+(V) = \left\{\lambda \in P(V)\,|\,\lambda + \eta \notin P(V) 
\ \mbox{and} \ \forall \eta >_0 0\right\}$.
\begin{lem}\label{lem2} Given $\lambda \in P(V)$, there exist $\eta \geq_0 0$ such that $\lambda +\eta \in P_+(V)$. 
In particular $P_+(V) \neq 0$.
\end{lem}
\begin{proof}
The proof proceeds as in \cite{E3}. 
For let $T= \{\mu \in P(V) \,|\, \mu (d)=\lambda (d)\}$.
Consider $W=\displaystyle{\bigoplus_{\mu \in T}} V_\mu $ which is $\gg(\overline{0},\overline{0})$ 
module and integrable. As 
in the proof of Lemma 2.6 of \cite{E3} one can prove $W$ is finite dimensional. Consider
$T_\lambda = \{\mu \in T \,|\, \lambda \leq_0 \mu \}$ which  is finite and hence 
has a maximal weight say $\gamma$. Now 
$\gamma - \lambda \geq 0$ and $(\gamma - \lambda) (d) =0.$ 
By observation \ref{ob1} we have $\gamma - \lambda \geq_0 0$. Take 
$\eta = \gamma-\lambda$. We need to prove $\lambda +\eta \in P_+(V)$. Suppose $\lambda+\eta+\eta_1 \in P(V)$ for 
$\eta_1 >_0 0$. But as $\eta_1(d_0)=0$,
it implies that $\lambda+\eta+\eta_1 \in T$ which contradicts the fact that $\lambda+\eta$ is maximal. 
This proves 
$\lambda+\eta \in P_+(V)$. Hence the lemma is proved.
\qed
\end{proof}
\begin{lem}\label{lem3}
Given $\lambda \in P_+(V),$ there exists $\lambda_1 = \lambda +P_1 m_0 
\delta + \eta_1 \in P(V)$ for $P_1 \in \NN$,
$\eta_1\geq 0$, and a $\lambda_1$ weight vector $u_1$ such that $N^+_{aff} u_1 =0$.
\end{lem}
\begin{proof}
The lemma follows from the proof of Theorem (2.4) of \cite{C}. The assumption of irreducibility 
is not needed for this part.
\qed
\end{proof}
\begin{prop}\label{pr3}
 Suppose for any $\lambda \in P(V)$, there exists $\eta >_1 0$ such that $\lambda+\eta \in P(V)$. 
Then there exists infinitely many $\lambda_i \in P(V), i \in \ZZ_+$ such that 
\begin{enumerate}
 \item $\lambda_i <_1 \lambda_{i+1},\lambda_{i+1} (d_0)-\lambda_i(d_0) \in \ZZ_+$ for $i \in \ZZ_+.$
 \item There exists $u_i \in V_{\lambda_i}$ such that $N^+_{aff} u_i =0$. 
 The corresponding irreducible module $V(\lambda_i) \subseteq V.$
 \item There exists a common weight for all $V(\lambda_i)$.
\end{enumerate}
\end{prop}
\begin{proof}
Choose $\lambda \in P_+(V)$. Then by Lemma \ref{lem3} there exists 
$\lambda_1 = \lambda + P_1 m_0 \delta +\eta_1 \in P(V)$ for $P_1 \in \NN$, $\eta_1 \geq 0$, 
and a weight vector $u_1$ of weight 
$\lambda_1$ such that $N^+_{aff} u_1 =0$. By assumption there exists $\stackrel{\sim}{\eta_1} >_1 0$ 
such that $\lambda_1+ \stackrel{\sim}{\eta_1} \in P(V)$. 
Now by Lemma \ref{lem2} there exists $\alpha \geq_0 0$ such that 
$\lambda_1+ \stackrel{\sim}{\eta_1}+\alpha  \in P_+(V)$. Using Lemma \ref{lem3} there exists 
$\lambda_2=\lambda_1+ \stackrel{\sim}{\eta_1}+\alpha+P_2 m_0 \delta +\eta_2 \in P(V)$
where $P_2 \in \NN$ and 
$\eta_2 \geq 0$. Further there exists $\lambda_2$ weight vector $\mu_2$ such that $N^+_{aff} u_2 =0$. 
Now by construction it is clear $\lambda_1 <_1 \lambda_2$. We claim that 
$\lambda_2(d_0) -\lambda_1 (d_0) > 0$. In any case 
$\lambda_2(d_0) -\lambda_1 (d_0) \geq 0$ as $A = \stackrel{\sim}{\eta_1}\!\!(d_0) 
+ \alpha(d_0) + P_2 m_0 + \eta_2(d_0) \geq 0$.
Suppose $A=0$. By using observation \ref{ob1} we see that  $B = \stackrel{\sim}{\eta_1} + \alpha + 
P_2 m_0 \delta + \eta_2 \in {\overset{\circ}{Q}^+}$. 
This contradict the fact that $\lambda_1 \in P_+(V)$ unless $B=0$. But $B \neq 0$ 
as $\stackrel{\sim}{\eta_1} >_1 0$.  
This proves our claim. By repeating this process infinitely many times we get (1) and (2).

Now we have $\lambda_i(d_0) = \lambda_1(d_0) + S_i$ where $\{S_i\}_, \,i \geq 2$ is 
strictly increasing sequence of positive integers. 
It is easy to see that there exists $j \in \NN,\, 0 \leq j < m_0$ such that $S_i \equiv j(m_0)$ for 
infinitely  many indices. Denote such 
sequence by $S_{i_g}, g =1,2,\ldots$ .
Write $S_{i_g} = f_{i_g} m_0 +j$ and can assume $f_{i_g} >0$. 
Let $\mu_0^j = \overline{\mu}_0 +(\lambda_1(d_0)+j) \delta +\lambda_1 (K) w$.
Consider
$$
\lambda_{i_g} - \mu_0^j = \overline{\lambda}_{i_g} - \overline{\mu}_0 + f_{i_g} m_0 \delta > 0.
$$
This proves $\mu_0^j \leq \lambda_{i_g}$ and hence $\mu_0^j \in P(V(\lambda_{i_g})$ by fact 2 given in 
the beginning of this subsection. For this subsequence 
(1), (2) and (3) are true. 
\qed
\end{proof}
{\bf Proof of Theorem \ref{thn2}.} \  Suppose $V$ is not a highest weight module. 
Then the assumptions of Proposition \ref{pr3} are true. 
Hence there exists infinitely many non-isomorphic irreducible 
modules with common weight. Then the dimension 
of the common weight is infinite dimensional contradicting the definition of weight module. 
Therefore $V$ is highest weight module. 
Let $\lambda$ be the top weight and put $N=V_\lambda$. Since $V$ is irreducible, 
by weight argument we see that $N$ is irreducible $\LL^0$- module. Thus $V \cong V(N)$.
\qed
\section{Reduction to non-graded module}\label{sec7}

 Recall the universal central extension of a Lie Tours $\stackrel{\sim}{LT}$ from 
 section \ref{sec2} and of $\LL$ from section \ref{sec3}. Let
$\stackrel{\sim}{LT}=\overline{LT}+\overline{D}$ where $\overline{D}$ is spanned
by derivations $d_0,d_1, \cdots, d_n$. Now 
let $D$ be the span of $d_1, \cdots, d_n$. In this section, for an irreducible 
weight module $V$ of $\stackrel{\sim}{LT}$, we associate an irreducible weight module 
$V_1$ for $\LL$ (need not be unique), 
 and in the next section we will show how to recover the original module $V$ from the
$\LL$ module $V_1$; there by proving that the 
study of irreducible weight modules for $\stackrel{\sim}{LT}$ is reduced to the 
irreducible weight modules for $\LL$. We 
generally refer $\LL$ modules as non-graded because the grade measuring $D$ is removed.
We fix an irreducible weight module $V$ of $\stackrel{\sim}{LT}$ for the rest of the section.
\subsection{}
\begin{defn} A linear map $z$ from $V\rightarrow V$ is called a central operator of 
degree $(k_0,k) \in \ZZ^{n+1}$ if 
$z$ commutes with $\overline{LT}$ and $d_i z-z d_i = k_i z$ for all $i$. 
\end{defn}
For example $t_0^{r_0} t^r K_j$ is a central operator 
of degree $(r_0,r)$ and $t^{s_0}_0 t^s K_i t_0^{r_0} t^r K_j$ is a 
central operator of degree $(r_0+s_0, r+s)$. The proof of the following lemma is exactly as 
given in Lemma 1.7 and Lemma 1.8 of \cite{E4}.
\begin{lem} \label{lem4}
 \begin{enumerate}\item Suppose $z$ is a central operator of degree $(k_0,k)$ 
 and suppose $z u \neq 0$ for some $u \in V$, then 
$z w \neq 0$ for all $w \in V$.
\item Suppose $z$ is a non-zero central operator of degree $(k_0,k)$, 
then there exists a non-zero central operators $T$ of degree 
$(-k_0,-k)$ such that $zT = Tz=Id$ on  $V$.
\item Suppose $z_1$ and $z_2$ are non-zero central operator of same degree,
then there exists a scalar $\lambda$ such that 
$z_1=\lambda z_2$.
\end{enumerate}
\end{lem}

  Let $S=\left\{(k_0,k) \in \ZZ^{n+1} \,|\, t_0^{k_0} t^k K_i \neq 0  \ \mbox{for some} 
  \ i\right\} \subseteq \Gamma \oplus \Gamma_0$.
Let $\langle S \rangle$ be the subgroup of $\Gamma \oplus \Gamma_0$ generated by $S$. 
Clearly rank $\langle S \rangle = p \leq n+1$. Note that for any 
$(k_0,k) \in \langle S \rangle$ there exists a nonzero central operator of degree $(k_0,k)$. 
Let $e_i = (0,\cdots,1,\cdots 0),  0 \leq i \leq n$ 
be the standard $\ZZ$- basis of $\ZZ^{n+1}$. Now it is standard fact that
upto a choice of co-ordinates, there exists 
positive integer $l_{n-p+1}, \cdots l_n, m_i|l_i$ such that 
$\{l_i e_i, i \geq n-p+1\}$ is a $\ZZ$- basis for $\langle S \rangle$.
\begin{prop}\label{props3} 
\begin{enumerate}
 \item There exists non-zero central operators $z_i, i \geq n-p+1$ of degree $l_i e_i$.
 \item $p < n+1$.
 \item Suppose $t_0^{k_0}t^k K_i \neq 0$ on $V$. Then $k_0=k_1= \cdots k_{n-p}=0$ and $i \leq n-p$.
 \item There exists $\overline{LT} \oplus D_p$ proper submodule $W$ of $V$ such that $V/W$ is a weight module for 
 $\overline{LT} + D_p$ with respect to $\hh(0) + \sum \CC K_i \oplus D_p$. Here $D_p$ is space spanned by 
$d_0,d_1, \cdots d_{n-p}$.
 \item We can choose $W$ such that $\{ z_i v-v \,|\, v \in V, \,n-p+1 \leq i \leq n \} \subseteq W$.
\end{enumerate}
\end{prop} 
The proof is exactly as given in Theorem (4.5) of \cite{E3}.

\begin{prop}\label{props2} Notations as in above; suppose $t_0^{k_0} t^k K_j$ acts non-trivially on 
$V$ for some $j$ and for 
some $(k_0,k) \in \ZZ^{n+1}$ then $p=n$.
\end{prop}
The proof is similar to Theorem 1.10 of \cite{E4}. We record the following which is of independent interest.
\begin{prop} Suppose $n \geq 1$, and $\stackrel{\sim}{V}$ is any weight module for $\stackrel{\sim}{LT}$ such that 
$t_0^k t^k K_i$ is zero on $\stackrel{\sim}{V}$ for $(k_0,k) \neq (0,0)$ and any $i$. Then each $K_i$ 
acts trivially on $V$.
\end{prop}
This follows from Corollary 1.24 \cite{E4}. Just note that $\stackrel{\sim}{V}$  is a module for 
the toroidal Lie algebra
$\gg(\overline{0},\overline{0}) \otimes \CC [t_0^{\pm m_0}, \cdots t_n^{\pm m_n}] \oplus z(m_0,m) 
\oplus {\overline{D}}$.
Also note that we are not assuming $\stackrel{\sim}{V}$  is irreducible.

 We have started with an irreducible weight module $V$ for $\stackrel{\sim}{LT}$. 
 Then we have proved it is actually 
a module for $\LL \oplus D = \, \stackrel{\sim}{\LL}$. See definition of $\varphi$ 
from beginning of section \ref{sec3} and
observe that $\ker \varphi$ is trivial on $V$ upto a choice of co-ordinates (see Prop. \ref{props3}(3) and
Prop. \ref{props2}).
Note that $D_n = \CC d_0$ and $V$ contains $\LL$ submodule $W$. We can assume that 
$\{ z_i v-v \,|\, v \in V, i\geq1 \} \subseteq W$.
\begin{prop}Notation as above; $V$ contains maximal $\LL$ submodule $W$ of $V$, so that $V/W$ is $\LL$- irreducible 
weight module. Further, we can choose $W$ so that $W$ contains $\{ z_i v-v \,|\, v \in V, i\geq1 \}$.
\end{prop}
The proof is similar to Theorem 2.5 of \cite{E4}.
\section{Passage from non-graded to graded}\label{sec8}  
The  main ideas comes from \cite{E4} but arguments differ.
 We will first record some general results on weight modules. 
Let $G$ be any Lie algebra and let $H$ be ad-diagonalizable 
subalgebra of $G$.
Let $G=\displaystyle{\bigoplus_{\alpha \in H^*}} G_\alpha$, where 
$G_\alpha= \{ x \in G \,|\, [h,x] = \alpha(h)x, \forall h \in H\}$.
Let $U$ be the universal enveloping algebra of $G$.
\begin{enumerate}
 \item Write $U=\displaystyle{\bigoplus_{\eta \in H^*}}U_\eta, U_\eta = 
 \{x \in U \,|\, [h,x]=\eta (h)x \, \forall\, h \in H \}$.
Clearly $U_0$ is a subalgebra of $U$ and contains $H$.
 \item Let $J$ be an irreducible weight module for $G$.
Write $J=\displaystyle{\bigoplus_{\lambda \in H^*}}J_\lambda, 
J_\lambda =\{ v \in J \,|\, h v =\lambda(h) v, \forall h \in H\}$.
Then it is standard fact that each $J_\lambda$ is an irreducible module for $U_0$. 
Fix $\lambda \in H^*$ such that $J_\lambda \neq 0$ and consider the induced module 
$$
M(\lambda) = U \displaystyle{\bigotimes_{U_0}} J_\lambda
.$$
\end{enumerate}
From the universal properties of induced modules, the following can be seen.

\subsection{}
\begin{prop} 
\begin{enumerate}
 \item $M(\lambda)_\lambda \cong J_\lambda$ as $U_0$- modules.
 \item $M(\lambda)$ has unique irreducible quotient $V(\lambda)$ and 
 $V(\lambda)_\lambda \cong J_\lambda$ as $U_0$- modules.
 \item For any $G$ weight modules $W$ such that $W_\lambda \cong J_\lambda$ as $U_0$- modules, 
 then $W$ is a quotient of 
$M(\lambda)$.
\end{enumerate}
\end{prop}
\begin{cor}\label{cor2} Suppose $J$ and $W$ are irreducible $G$ weight 
modules such that $J_\lambda \cong W_\lambda$ as 
$U_0$ modules. Then $J \cong W$ as $G$-modules.
\end{cor}

  We will now continue with section \ref{sec7}. We have an irreducible weight module 
  $V$ for $\stackrel{\sim}{LT}$. We have 
non-zero central operators $z_1, \cdots, z_n$ of degree $l_i e_i, m_i | l_i$. We also have 
a maximal $\LL$ submodule $W$ of $V$ and $V_1=V/W$ is an irreducible weight module. 
We assume $W$ contains the space spanned by 
$\{z_i v-v \,|\, v \in V, 1 \leq i \leq n\}$. Also we have the group 
$\langle S \rangle = \ZZ l_1 \oplus \cdots \oplus \ZZ l_n \subseteq \Gamma$. For 
every $s\in\langle S \rangle$, there is a non-zero central operator of degree $s$. Let 
$[S]$ be the space of such operators. As operators 
on $V, K \otimes A(m) \subseteq [S]$. Equality may not hold. Let 
$\overline{G} = \overline{LT} +[S] \oplus \CC d_0$ which is a Lie algebra 
in an obvious sense and contains $\LL$.
Let $G= \overline{G} \oplus D$ and contains $\stackrel{\sim}{LT}$. 
Let $H= \hh(0) \oplus \sum \CC K_i + \CC d_0 +D = \,\stackrel{\sim}{\hh},$ be an ad-diagonalizable
subalgebra of $G$. We note that $V$ is 
irreducible $G$ module and $V_1=V/W$ is an irreducible 
$\overline{G}$-module. Write $\overline{U}= U(\overline{G})$ and $U=U(G)$ 
which are universal enveloping algebra of 
$\overline{G}$ and $G$. 
Note that by our assumption, each $z_i$ acts as 1 on $V_1$. It is easy to 
see that $z_i^{-1}$ also acts as 1 on $V_1$. Suppose 
$z$ is any non-zero central operator of degree $k \in \langle S \rangle$. Then $ z = \lambda 
\Pi z_i^{p_i}, p_i \in \ZZ, \lambda \neq 0$ and 
$k=(p_1, \cdots, p_n)$. This can be seen from Lemma \ref{lem4}(3).

 Let $\LL(V_1) = V_1 \otimes A_n$. Recall that $\overline{G}$  is 
 $\ZZ^n$- graded via $D$; 
write $\overline{G} =\displaystyle{\bigoplus_{k \in \ZZ^n}} \overline{G}_k$.
Define $G$-module action on $(\LL(V_1), \Pi (\alpha))$ where 
$\alpha =(\alpha_1, \cdots, \alpha_n) \in \CC^n$ by
$$
\begin{array}{lll}
 X.v(r) &=& (Xv) (r +k), X \in  \overline{G}_k;\\
 d_0.(v(r)) &=& (d_0 v)(r);\\
 d_i. v(r) &=& (r_i+\alpha_i) v (r), 1 \leq i \leq n, r \in \ZZ^n.
\end{array}
$$
It is easy to see that $\LL(V_1)$ is a  $G$-module and each central operator from $[S]$ acts 
as invertible operator. 
Any $\lambda \in \, \stackrel{\sim}{\hh}\!\!(0)^*$ can be extended to 
$\stackrel{\sim}{\hh}$ by defining  $\lambda(D)=0$ and $\lambda(K_i) = 0$ for $1 \leq i \leq n$. 
Recall the definition 
$\delta_r \in \, \stackrel{\sim}{\hh}\!\!(0)^*$ from section \ref{sec2}.
Clearly $\LL(V_1)_{\lambda+\delta_r}=(V_1)_\lambda \otimes t^r$
which is finite dimensional. Thus $\LL(V_1)$ is a $G$ weight module.
In this section our aim is to prove that $\LL(V_1)$ is direct sum of finitely many 
irreducible $G$-modules. All components 
are isomorphic upto grade shift. Further we will prove that one of the
components is isomorphic to $V$ as $G$-modules for 
suitable $\alpha$. There by recovering the original module $V$ from $V_1$.
 Consider $\rho:\LL(V_1) \rightarrow V_1$ defined by $\rho(v 
\otimes t^r) = v$ which is a surjective $\overline{G}$- module map. 
Let $v(r) = v \otimes t^r$.
\begin{lem}\label{lem5} Suppose $W_1$ is a non-zero submodule of $\LL(V_1)$. Then 
\begin{enumerate} 
\item $\rho(W_1) =V_1$.
\item  Suppose $\lambda \in P(V_1)$ then $\lambda + \delta_k \in P(W_1)$ for some $k \in \ZZ^n$.
\end{enumerate}
\end{lem}
\begin{proof} Since $V_1$ is irreducible, to see (1) it is sufficient to prove that $\rho(W_1) \neq 0$. But $W_1$ is a 
$D$- module and hence contains vectors of the form $v(k)$. Then $\rho(v(k) = v \neq 0$.
To see (2), let $v \in V_1$ of weight $\lambda \in \stackrel{\sim}{\hh}\!\!(0)^*$. Then there exists $w \in W_1$ such that 
$\varphi(w) = v$ by (1). Write $w= \sum w_i (\gamma^i)$ for some $w_i \in V_1$ and $\gamma^i \in \ZZ^n$. Each $w_i$ is 
weight $\lambda$. As $w_1$ is $D$- module, $w_i (\gamma^i) \in W_1 \ \mbox{for all} \ i$. 
Thus $\lambda + \delta_{\gamma i} \in P(W_1)$.
\qed
\end{proof}
{\bf Remark:}
 Let $W^1$ be a non-zero submodule of $\LL(V_1)$. Then for any $z \in [S]$ of degree $k$, we have 
$z W^1_\lambda = W^1_{\lambda +\delta_k}$. This is because $W^1$ is $[S]$ invariant and $z$ is invertible.

Fix a $\lambda \in P(V_1)$ and define $K(W^1) = \displaystyle{\sum_{0 \leq s_i < l_i}} \dim 
(W^1_{\lambda +\delta_s})$ for any submodule $W^1$ of $\LL(V_1)$. We claim that $K(W^1) >0$ for non-zero submodule 
$W^1$. From Lemma \ref{lem5}(2), $\lambda + \delta_k \in P (W^1)$ for some $k \in \ZZ^n$. Now by above remark  
we can assume $0 \leq k_i < l_i$. This proves $K(W^1) >0$. 
Suppose $W^1 \subseteq W^2$ are submodules then clearly 
$K(W^1) \leq K(W^2)$. 
\begin{lem}\label{lem6} Let $W^1$ be a non-zero $G$ submodule of $\LL(V_1)$. Then $W^1$ contains an irreducible 
$G$- submodule.
\end{lem}
\begin{proof} All submodules consider here are non-zero $G$-modules.  
Choose a submodule $W^2$ of $W^1$ such that $K(W^2)$ is minimal .
Let $Q=\displaystyle{\bigoplus_{0 \leq k_i < l_i}} W^2_{\lambda +\delta_k}$ and 
let $\stackrel{\sim}{W}$ is a submodule 
generated by $Q$.
Clearly $K(\stackrel{\sim}{W}) = K(W^2)$. We claim that $\stackrel{\sim}{W}$ is irreducible.
For that let $W^3$ a submodule 
of $\stackrel{\sim}{W}$. We have 
$\displaystyle{\bigoplus_{0 \leq k_i < l_i}} W^3_{\lambda+\delta_k} 
\subseteq \displaystyle{\bigoplus_{0 \leq k_i < l_i}} 
W^2_{\lambda + \delta_k}$. By minimality we know $K(W^3) = K(W^2)$ and 
hence equality holds above. Thus $W^3$ contains $Q$ and hence 
contains $\stackrel{\sim}{W}$. Thus $W^3 = \stackrel{\sim}{W}$. 
This proves that  $\stackrel{\sim}{W}$ is irreducible.
\qed
\end{proof}
{\bf Remark:}
For $v \in V_1$ and $k \in \ZZ^n$ let $U v(k)$ be the $G$ submodule of $\LL(V_1)$ generated by $v(k)$.
 Let $v$ be a weight vector of $V_1$ and let $r$ and $s \in \ZZ^n$.
Consider the map
$
\varphi : Uv(r) \rightarrow Uv(s) \ \mbox{by} \ \varphi w(k)=w(k+s-r)
$
for $w(k) \in U v(r)$ where $w \in V_1$ and $k \in \ZZ^n$. It is easy to check that $\varphi$ 
is $\overline{G}$ module 
isomorphism. But need not be $G$-module map. Suppose $D$ acts on $v(r)$ by $\alpha \in \CC^n$. 
That means $d_i v (r) = \alpha_i v(r)$. Suppose $D$ acts on $v(s)$ by $\beta$. 
Then it is easy to see that 
$U v(r) \cong Uv(s) \otimes \CC$ as $G$-modules where $\CC$ 
is a one dimensional module for $G$ with $\overline{G}$ acting 
trivially on $\CC$ and $D$ acts as $\alpha - \beta$. 
This is what we call $Uv(r)$ isomorphic to $Uv(s)$ as $G$-module upto a grade shift. 
Now it is easy to see that $Uv(r)$ is irreducible if $Uv(s)$ is irreducible as $G$-module.

\begin{prop} There exists a weight vector $v$ in $V_1$ such that 
\begin{enumerate}
 \item $U v(k)$ is an irreducible $G$-module for all $k \in \ZZ^n.$
 \item $\displaystyle{\sum_{k \in \ZZ^n}}Uv(k) = \LL(V_1).$
 \item $\displaystyle{\sum_{0 \leq k_i < l_i}}Uv(k) = \LL(V_1)$.
\end{enumerate}
\end{prop}
\begin{proof}
  Note that (3) follows from (2) as $Uv(k) = Uv(k+r)$ for any 
$r \in \langle S \rangle$. One can use central operators of degree $r$. See remark followed by Lemma \ref{lem5}.
\begin{enumerate}
 \item Let $W^1$ be an irreducible  $G$-module of $\LL(V_1)$ (Lemma \ref{lem6}). 
 From the proof of Lemma \ref{lem5}(1) we see that $W^1$ 
contains a vector $v(k)$ for some $k \in \ZZ^n$, and some weight vector $v$ in $V_1$. 
Now $Uv(k) \subseteq W^1$ and hence 
$Uv(k)$ is irreducible and $Uv(k) =  W^1$. Now by above remark it follows that 
$Uv(s)$ is irreducible for any $s \in \ZZ^n$.
\item Let $w(s) \in \LL(V_1)$ for $w \in V_1$ and $s \in \ZZ^n$. 
Since $V_1$ is $\overline{G}$ irreducible and $v \in V_1$ 
there exists $X \in \overline{U}$ such that $X v =w$.
Write $X = \sum X_{k^i}$ where $D$ acts on $X_{k^i}$ by $k^i$.
Consider $\sum \pi(\alpha) X_{k^i} v (s-k_i) = \sum (X_{k_i} v) (s) = w(s)$ 
(recall the definition of $G$ action on 
$\LL(V_1)).$  Thus $\LL(V_1)\subseteq \displaystyle{\sum_{s \in \ZZ^n}} U v(s)$ and hence we have (2). 
\end{enumerate}
\qed
\end{proof}
\begin{thm}
\begin{enumerate}
 \item $(\LL(V_1), \pi(\alpha))$ is completely reducible as $G$-modules and the number of components 
are finite. Further all components are $G$- isomorphic upto grade shift. 
 \item For a suitable $\alpha, a$ component of $(\LL(V_1), \pi(\alpha))$ is isomorphic to $V$ as $G$-modules.
\end{enumerate}
\end{thm}
\begin{proof} We have already seen that $\displaystyle{\sum_{0 \leq s_i < l_i}} Uv(s) = \LL(V_1)$  
for some weight vector $v$ in $V_1$. Let $T=\{ s \in \ZZ^n \,|\, 0 \leq s_i < l_i \}$.
Suppose for some $\gamma$,
$Uv(\gamma) \cap \displaystyle{\sum_{\substack{ s \in T \\ s \neq \gamma}}}
U v(s) \neq 0$. Since $U v(\gamma)$ is irreducible 
it follows that $Uv(\gamma) \subseteq \displaystyle{\sum_{\substack{ s \in T \\ 
s \neq \gamma}}} U v(s) \neq 0$. Thus we have 
reduced one component from the sum. By repeating this 
finitely many times we see that $\LL(V_1)$ is direct sum of fewer terms. 
This proves 1.
To prove 2) recall that $U$ and $\overline{U}$ are universal enveloping algebra 
of $G$ and $\overline{G}$ respectively. 
Recall the notation from beginning of this section. We have $U_0= \overline{U}_0. U(D)$. 
Consider the weight space 
$V_{\lambda+\delta_\gamma} \subseteq V$ which is an irreducible $G$-module. 
As noted earlier $V_{\lambda+\delta_\gamma}$ is an 
irreducible $U_0$- module and remains irreducible as $\overline{U}_0$ 
as $U(D)$ acts as scalars. Recall the map 
$\rho :\LL(V_1) \rightarrow V_1$ defined by $\rho(v(r)) = v$. Now $\rho$ restricted to 
$V_{\lambda+\delta_\gamma}$ is injective as $\ker \rho$ 
cannot contain $\stackrel{\sim}{\hh}$ weight vectors.
Now $\rho(V_{\lambda+\delta_\gamma}) \subseteq (V_1)_\lambda$ is an 
irreducible $\overline{U}_0$ module and $\rho$ is a 
$\overline{U}_0$ module map. Therefore $\rho(V_{\lambda+\delta_{\gamma}}) = (V_1)_\lambda$
and hence $\rho(V_{\lambda+\delta_{\gamma}}) \cong (V_1)_\lambda$ as $\overline{U}_0$- modules.
As $V_1 \subseteq \LL(V_1)$ we see that $(V_1)_\lambda = \LL(V_1)_\lambda 
= \oplus (Uv(\gamma))_\lambda$ as $\overline{U}_0$- modules. 
Each $(Uv(\gamma))_\lambda$ is $U_0$ irreducible and remains 
irreducible as $\overline{U}_0$- modules. This proves 
$V_{\lambda+\delta_\gamma} \cong (Uv(s))_\lambda$ as 
$\overline{U}_0$- modules for some $s \in T$.
We know that $D$ acts as scalars on $V_{\lambda+\delta_\gamma}$ and 
$(Uv(s))_\lambda$. We can choose $\alpha$ in such a way that
$D$ action is same on $V_{\lambda+\delta_s}$ and $(Uv(s))_\lambda$. 
Thus for the special choice of $\alpha$ we see that 
$V_{\lambda+\delta_s} \cong (Uv(s))_\lambda$ as $U_0$- modules. 
Now by Corollary \ref{cor2} it follows that 
$V \cong Uv(s)$ as $G$-modules. This completes the proof of (2).
\qed
\end{proof}
{\bf Remark:} We would like to mention an important result Theorem 5.14 in 
\cite{NK}. It is proved that a multi-loop algebra 
with  $\gg(\overline{0},\overline{0}) \neq 0$ is support-isomorphic to a Lie torus. 
It is easy to see that real roots goes to 
real roots under support-isomorphism. Hence the integrability carries on.
It is not clear to the author how to extend 
support-isomorphism to their universal central extension. 
If such a lift exists, then all our classification results holds 
for the multi-loop algebra such that  $\gg(\overline{0},\overline{0}) \neq 0$.

\bibliography{Integrablemodules.bib}

\nocite*{}

\bibliographystyle{plain}
\end{document}